\documentclass[a4paper]{article}
\usepackage{theorem}
\usepackage{amssymb}
\usepackage{times}
\usepackage{mathptm}
\title{Oriented Lagrangian Matroids}
\date{}
\author{Richard F. Booth\footnote{\texttt{\{richard.booth\}\{alexandre.borovik\}@umist.ac.uk}}
\and  Alexandre V. Borovik\footnotemark[1] 
\and  Israel M. Gelfand\footnote{{\tt igelfand@math.rutgers.edu}}
\and Neil White\footnote{{\tt white@math.ufl.edu}}
\,\thanks{Partially supported by EPSRC Grant GR/M24707.}}
\newcommand{\qed}{\hfill $\diamond$ \\ \vskip 5pt}
\theorembodyfont{\slshape}
\newtheorem{lem}{Lemma}
\theorembodyfont{\slshape}
\newtheorem{thm}[lem]{Theorem}
\newtheorem{cor}[lem]{Corollary}
\theorembodyfont{\upshape}
\newtheorem{defn}{Definition}
\theorembodyfont{\slshape}

\theorembodyfont{\slshape}
\newtheorem{ax}{Axiom}
\newcommand{\lew}{\leqslant^w}
\newcommand{\B}{\ensuremath{{\cal B}}}
\newcommand{\K}{\ensuremath{\mathbb{K}}}
\newcommand{\R}{\ensuremath{\mathbb{R}}}
\newenvironment{proof}[1][Proof]{\paragraph{#1}}{\qed}
\begin{document}

\pagestyle{myheadings}
\markright{{\small \sc R.~F.~Booth et al. \qquad \qquad \qquad \qquad
Oriented Lagrangian Matroids }}
\maketitle

The aim of this paper is to develop, by analogy with the theory of
oriented (ordinary) matroids \cite{oribook}, an oriented version of
the theory of Lagrangian symplectic matroids \cite{bgw}. Recall that
the concept of an oriented matroid axiomatises, in combinatorial form,
the properties of the ensemble of non-zero $k \times k$ minors of a
real $k \times n$ matrix $A$, say. The $k$-subsets of indices of
columns of $A$ such that the corresponding $k \times k$ minors are
non-zero form the collection of bases of an (ordinary) matroid
\cite{whitebook}. When we consider also the signs of these minors we
obtain a considerably richer combinatorial structure, an (ordinary)
oriented matroid. These structures have important applications in
combinatorics, geometry and topology. To give an idea of the concept
of an oriented Lagrangian matroid, we consider one of the simplest
situations in which they arise.

Assume that we have a symmetric $n \times n$ matrix $A$ over some
field \K. The collection of sets of row (equivalently, column)
indices $K \subseteq \{\,1,\ldots, n\,\}$ corresponding to
non-zero diagonal minors $|a_{ij}|_{i,j \in K}$ of $A$  forms what
is known as a $\Delta$-matroid \cite{bou3} (we include for formal
reasons the empty set $K = \emptyset$, and specify that a $0
\times 0$ minor takes the value $1$). These are equivalent
structures to Lagrangian symplectic matroids. When $\K=\R$, we may
consider also the signs of these minors. What follows in this
introduction is a very elementary description of the resulting
combinatorial structure; all necessary proofs, and detailed
explanations, may be found later in the paper.

Let our symmetric $n \times n$ matrix $A$ over \K\  have columns
indexed by $I=1, \ldots, n$, and let $e_1, \ldots, e_n$ denote the
standard basis in $\R^n$.  Given some $B \subseteq I$ we define a
point in $\R^n$ by
$$ B \rightarrow \sum_{i \in B} e_i - \sum_{i \notin B} e_i, $$ a
vertex of the $n$-dimensional hypercube $H=[-1, 1]^n$. Consider the
convex hull of the points induced from non-zero minors of $A$ in this
way. It can be shown that this polytope, $P$, satisfies:
\begin{description}
\item[(*)] All vertices of $P$ are vertices of $H$, and edges of $P$ are
either
edges of $H$ or diagonals of 2-dimensional faces of $H$.
\end{description}
Notice that if points corresponding to $B,C \subseteq I$ are at
opposite ends of an edge of $H$, without loss of generality $C = B
\sqcup \{i\}$ for some $1 \leqslant i \leqslant n$. We now
introduce some extra structure when $\K=\R$ by orienting edges of
$P$ which coincide with those of $H$. Given $B, C$ as above, we
direct the corresponding edge  from $B$ to $C$ if the
corresponding minors have the same sign, and  from $C$ to $B$ if
the signs differ. It will be shown in Theorem~\ref{thm:oripol}
that the resulting partially oriented $2$-skeleton of the polytope
satisfies
\begin{description}
\item[$(**)$] In any 2-dimensional face, the same number of edges are
oriented
clockwise as anti-clockwise.
\end{description}
We may abstract this situation as follows. $(*)$ is a
cryptomorphic definition of a Lagrangian symplectic matroid
\cite{bgw}, and $(**)$ is one of a pair of equivalent definitions
of orientations of  oriented Lagrangian matroids given in this
paper. This is a natural definition, as the following observation
shows: when the construction from a symmetric matrix is followed
the signs of all diagonal minors may be recovered from the
partially oriented polytope. Consequently, the rank and signature
of the corresponding quadratic form are preserved within the
combinatorial information (Corollary \ref{cor:last}). In both
cases, the recovery of this information is both geometric and
uncomplicated.

The reader should not be surprised by the fact that we define new
combinatorial concepts in geometric terms. We work in the
framework of the theory of Coxeter matroids, of which (ordinary)
matroids and Lagrangian symplectic matroids are special cases
corresponding to the Coxeter groups $A_n$ and $BC_n$.
Combinatorial properties of Coxeter matroids are easily translated
into the properties of the corresponding polytopes, and vice versa
\cite{polytopes,borovik-vince}. In particular, the adjacency of
vertices of the polytope $P$ can be expressed entirely in terms of
the Lagrangian symplectic matroid associated with the matrix $A$
(Theorem~\ref{thm:adjcrit} below; it is a special case of
\cite{borovik-vince}). Moreover, the combinatorial type of the
polytope $P$ (in particular, its faces) can be read from the
corresponding Lagrangian symplectic matroid \cite{book}.

One of the uses of oriented matroids is to provide a finer
stratification of the real Grassmannian $G_{n, k}(\R)$ than that
given by thin Schubert cells
\cite{bor-gel:schu,gelfand-serganova}. If $X$ and $Y$ are two
points in $G_{n, k}(\R)$ in the same thin Schubert cell $C$ (which
is equivalent to saying that they represent the same matroid), but
which represent different oriented matroids, then they belong to
different connected components of $C$. This finer stratification
 was used by MacPherson in his construction of
the combinatorial model for the  Grassmannian \cite{macpherson}.

Similarly, oriented Lagrangian symplectic matroids provide a finer
stratification of the variety of maximal isotropic subspaces of a
real symplectic $2n$-dimensional vector space (see Sections 1.2
and 3). Points in a thin Schubert cell representing different
oriented Lagrangian symplectic matroids lie in distinct connected
components of the thin Schubert cell. It would be interesting to
use our concept of oriented Lagrangian symplectic matroid for the
construction, by analogy with \cite{macpherson}, of a
combinatorial model for the variety of maximal isotropic subspaces
of a real symplectic $2n$-dimensional vector space.

We wish to emphasise that oriented Lagrangian symplectic matroids
should not be viewed as generalisations of oriented ordinary matroids.
In fact, the natural generalisation of oriented matroids is provided
by oriented even $\Delta$-matroids \cite{wenz,w3} (or oriented
Lagrangian orthogonal matroids, in terminology of \cite{boo2}). Every
(ordinary) matroid is a Lagrangian orthogonal matroid, and every
Lagrangian orthogonal matroid is also a Lagrangian symplectic matroid,
which can be explained from the embeddings of root systems $A_{n-1} <
D_n < C_n$ \cite{bgw,book}. However, the concepts of orientation for
Lagrangian orthogonal and symplectic matroids are very different
because they reflect the geometric properties of two different
geometric objects: the varieties of maximal isotropic subspaces in the
space $\R^{2n}$ with a non-degenerate symmetric or a skew symmetric
form. Not surprisingly, (ordinary) oriented matroids fit more happily
into the geometry of a symmetric scalar product. See \cite{boo2} for
more detail on orthogonal orientation.

\section{Symplectic Matroids}

This exposition of symplectic matroids follows \cite{bgw}. For a more
general
discussion of the theory of Coxeter
matroids, of which symplectic matroids are a special case, see
\cite{book}.

\subsection{Symplectic and Lagrangian Matroids}

Let
$$ I=\{1,2, \ldots, n\}, \quad I^*=\{1^*,2^*,\ldots,n^*\} ,$$ and $J=I
\cup I^*$. We define a pair of maps $^*:I \rightarrow I^*$ by $ i
\rightarrow i^*$ and $^*:I^* \rightarrow I$ by $ i^* \rightarrow i$,
giving us $^*: J \rightarrow J$ an involutive permutation of $J$. We
say that a subset of $K \subset J$ is \emph{admissible} if and only if
$K \cap K^* = \emptyset$. We denote by $J_k$ the collection of all
admissible $k$-subsets of $J$. (We often use the word `collection' to
denote a set of sets, in order to avoid confusion.)

An admissible permutation $w$ of $J$ is one which satisfies
$w(i^*)=w(i)^*$ for each $i\in J$. We call the group of all admissible
permutations $W$, and remark that it is the \emph{hyperoctahedral
group} $BC_n$, that is, it is isomorphic to the group of symmetries of
the $n$-cube $[-1,1]^n$ in the real Euclidean space $\R^n$. Consider:
with $e_1, \ldots, e_n$ the standard orthogonal basis in $\R^n$, we
see that $W$ acts as follows: for $i \in I$ we set $ e_{i^*} = -e_i$
and $ we_i = e_{w(i)}$. It is not hard to see that $W$ is exactly the
group of all orthogonal transformations of $\R^n$ preserving the set
of vectors $\{ \pm e_1, \ldots, \pm e_n \}$ and thus the $n$-cube.

We order $J$ by $$ n^* < (n-1)^* < \ldots < 1^* < 1 < 2 < \ldots <
n $$ and for each $w \in W$ define a new ordering $\lew$ on $J$ by
putting $$ i \lew j \mbox{ if and only if } w^{-1} i \leqslant \
w^{-1} j. $$ An admissible ordering $\prec$ of $J$ is defined by $
i \prec j$ implies $j^* \prec i^*$, so we see that for $w \in W$,
$\lew$ is admissible. Now extending the concept of orderings to
elements of $J_k$, we take $\prec$ an admissible ordering on $J$,
and let $A,B \in J_k$ with $$ A= \{ a_1 \prec a_2 \prec \ldots
\prec a_k\} \mbox{ and } B=\{b_1 \prec b_2 \prec \ldots \prec
b_k\}. $$ Then we set $A \prec B$ if and only if $a_i \prec b_i$
for each $i=1, \ldots, k$. Take $\B \subseteq J_k$ and let $M$ be
the triple $(J,^*,\B)$. Then $M$ is a symplectic matroid if it
satisfies Axiom \ref{maxprop} below:

\begin{ax}[Maximality Property] \label{maxprop}
For every $w \in W$ there exists $A \in \B$ such that $B \lew A$ for every
$B
\in \B$.
\end{ax}
i.e. there is a unique $w$-maximal element of \B\ for each $w\in W$. If
the
axiom holds, we call \B\ the collection of \emph{bases}
of $M$. The cardinality $k$ of the bases is called the \emph{rank} of
$M$. A
symplectic matroid may be considered as a
Coxeter matroid for $BC_n$.

A \emph{Lagrangian} matroid is a symplectic matroid of rank $n$; this
means
that for each $i \in I$, we have either
$i \in A$ or $i^* \in A$ for every $A \in \B$. This paper is concerned
primarily with Lagrangian matroids.

\subsection{Representability}

Let $A,B$ be a pair of $k\times n$ matrices over a field $F$. Let $C$
be the $k\times 2n$ matrix $(A,B)$. Assume $C$ has rank $k$. We label the
columns of $A$ with
$I$ in order, and similarly those of $B$ with $I^*$, so that the
columns of $C$ are indexed by $J$. Define $\B \subseteq J_k$ by $X \in \B$
if and only if
\begin{enumerate}
  \item $X \in J_k$ and
  \item the $k\times k$ minor consisting of the $j$-th column of $C$ for
all $j
\in X$ is non-zero.
\end{enumerate}
Then we have

\begin{thm}
If $AB^t$ is symmetric, then \B\ is the collection of bases of a
symplectic
matroid.
\end{thm}

This result is proven in \cite{bgw}. We call a symplectic matrix arising
from a
matrix $(A,B)$ with $AB^t$ symmetric a
\emph{representable} symplectic matroid, and say that $(A,B)$ is a
representation of it over $F$. We sometimes refer to
this as a symplectic representation, to distinguish it from other kinds of
representations. Also, note that a representable
symplectic matroid is Lagrangian if and only if $A$ and $B$ are both $n
\times
n$ square matrices and $(A, B)$ is full rank.

Note that this result may be given in terms of isotropic subspaces of
symplectic vector spaces.
Let $V$ be a vector space with basis $e_1, \ldots,
e_n, {e_1}^*, \ldots, {e_n}^*$. Define the standard symplectic bilinear
form on
$V$ by putting
$$<e_i,e_j>=\left\{\begin{array}{lr}
1 &\hbox{ if $i\in I$ and $j=i^*$.}\\
-1 &\hbox{ if $j\in I$ and $i=j^*$.}\\
0 &\hbox{ otherwise.}
\end{array}\right.$$
Consider an isotropic subspace $U$ of $V$ (that is, $<u,v>=0$ for all
$u,v\in
U$). Take a basis of $U$ and write these as
row vectors in terms of the above basis in the above order,
thus forming a $k\times 2n$ matrix $C$ of rank $k$. This is
then a representation of a symplectic matroid, as above; the isotropic
condition is equivalent to the symmetry of $AB^t$.
Simple matrix algebra shows that the matroid represented is dependent only
on
the subspace chosen; thus, we may speak
of a subspace being a representation. A Lagrangian subspace (maximal
isotropic)
corresponds to a Lagrangian matroid.

We now define the concept of a thin Schubert cell.
Consider the variety of all Lagrangian subspaces of $V$.
Each such subspace has a corresponding Lagrangian matroid,
and the set of subspaces corresponding to the same matroid
is called a \emph{thin Schubert cell}.
Thin Schubert cells over $\mathbb{R}$ are in general not connected in the
real
topology.
Indeed, two subspaces $U$ and $U'$ corresponding to the same symplectic
matroid,
but with distinct signs of non-zero $k \times k$ minors
in their matrices $C$ and $C'$ cannot be continuously transformed one
into another within the same thin Schubert cell and thus belong to
distinct  connected components.

This concept of representation is related to one of several given in
\cite{bou3}, in which a representation consists of a symmetric,
$n\times n$ matrix A, and is equivalent to our concept of a representation
with
the additional requirement that $B=I_n$,
the identity matrix,  and that $k=n$. Our Lagrangian matroids can be shown
to
be equivalent to constructions referred to
as \emph{symmetric} matroids in \cite{bou3}; see \cite{wenz} for a proof
of
this.

\subsection{Symplectic Matroid Polytopes}

We now return to the $n$-cube $[-1,1]^n$ in $\R^n$. We again set $e_1,
\ldots,
e_n$ to be the standard orthonormal basis
in $\R^n$ and define $e_{1^*}=-e_1, \ldots, e_{n^*} = -e_n$. Now given $A
\in
J_k$ we set
$$e_A = \sum_{j \in A} e_j \in \R^n$$
so that, for example, the set $\{e_A | A \in J_n \}$ is the set of
vertices of
the $n$-cube.
We define the \emph{matroid polytope} of a symplectic
matroid $M$ to be the convex hull of the points $e_A, A \in \B$.
For example, the matroid polytope corresponding to the matroid represented
by
the matrix
$$ \left( \begin{array}{ccc|ccc}
        1 & 1 & 1 & 1 & 0 & 0 \\
        1 & 2 & 2 & 0 & 1 & 0 \\
        1 & 2 & 2 & 0 & 0 & 1
\end{array} \right) $$
can be drawn as shown in Figure 1,
where dots represent bases, solid lines edges of the matroid polytope and
dotted lines edges of the $n$-cube not present in the polytope.
The bases of this matroid are
$$ 1^*2^*3^*,12^*3^*,1^*23^*,1^*2^*3,123^*,12^*3 $$
(notice that we write, for example, $123$ rather than $\{1,2,3\}$.
This abbreviated notation is usual in matroid theory).
We shall be specifically interested in Lagrangian matroid polytopes,
that is, matroid polytopes corresponding to Lagrangian matroids.
\begin{figure} \caption{}

\unitlength 1mm
\begin{picture}(80,80)(-20,0)
{\tiny
\put(9,11){\makebox(0,0)[br]{$(1,-1,-1)$}}
\put(51,9){\makebox(0,0)[tl]{$(1,1,-1)$}}
\put(9,51){\makebox(0,0)[br]{$(1,-1,1)$}}
\put(51,51){\makebox(0,0)[tl]{$(1,1,1)$}}
\put(29.5,25.5){\makebox(0,0)[br]{$(-1,-1,-1)$}}
\put(29,66){\makebox(0,0)[br]{$(-1,-1,1)$}}
\put(71,21){$(-1,1,-1)$}
\put(71,66){$(-1,1,1)$}
}
\thinlines
\qbezier[20](10,50)(30,50)(50,50)
\qbezier[20](50,10)(50,30)(50,50)
\qbezier[20](30,65)(50,65)(70,65)
\qbezier[20](70,25)(70,45)(70,65)
\qbezier[15](50,50)(60,57.5)(70,65)
\thicklines
\put(10,10){\line(0,1){40}}
\put(10,10){\line(1,0){40}}
\put(30,25){\line(0,1){4}}\put(30,31){\line(0,1){34}}
\put(30,25){\line(1,0){4}}\put(36,25){\line(1,0){34}}
\put(10,10){\line(4,3){20}}                     
\put(50,10){\line(4,3){20}}
\put(10,50){\line(4,3){20}}
\put(10,50){\line(1,-1){40}}
\put(30,65){\line(1,-1){40}}
\put(10,10){\circle*{1.49}}
\put(50,10){\circle*{1.49}}
\put(10,50){\circle*{1.49}}
\put(30,25){\circle*{1.49}}
\put(30,65){\circle*{1.49}}
\put(70,25){\circle*{1.49}}
\thinlines
\end{picture}

\end{figure}

Some of the properties of these structures follow. The points
$e_A$ for $A \in \B$ are exactly the vertices of the matroid
polytope. In the sequel we shall often simply write $A$ in the
abbreviated notation used above rather than $e_A$. The matroid
polytope is inscribed in the $n$-cube, and in the Lagrangian case
its vertices are vertices of the $n$-cube also. The edges of a
Lagrangian matroid polytope can be of just two types, which may be
distinguished by their lengths.
 Type 1 is of length 2, for example the edge from $12^*3$ to
${1^*2^*3}$ in Figure 1. If $A$ and $B$ the bases corresponding to
the end vertices of an edge of type 1, then they are connected by
a {\em short exchange}, $A = B \Delta \{i,i^*\}$ for some $i$.
Type 2 is of length $2\sqrt 2$, for example the edge from $12^*3$
to $123^*$ in Figure 1. In that case, the end bases $A$ and $B$ of
the edge are connected by a {\em long exchange}: $A = B
\Delta\{i,j,i^*,j^*\}$ with $i\ne j,j^*$.  In fact, a polytope is
a Lagrangian matroid polytope in this sense if and only if its
vertices are vertices of the $n$-cube and its edges are only of
lengths 2 and $2 \sqrt 2$. For more details and proofs see
\cite{bgw,car}. These results rely heavily on the
Gelfand--Serganova theorem, of which the relevant special case for
symplectic matroids is Theorem 10 in \cite{bgw}. We shall call an
edge of type 1 a short exchange or short edge and an edge of type
2 a long exchange. Furthermore, it is easily seen that every edge
in the Lagrangian case lies upon the surface of the $n$-cube, that
is within one of the 2-dimensional faces of the $n$-cube.  In
fact, it is not hard to show
\begin{thm}
\label{thm:adjcrit}
Two bases $A,B$ of a Lagrangian matroid are adjacent if and only if
\begin{enumerate}
\item They are connected by a short exchange ($A=B\Delta\{i,i^*\}$);
or
\item They are connected by a long exchange
($A=B\Delta\{i,j,i^*,j^*\}$ with $i\ne j,j^*$) and
\[ |\{ B\Delta\{i,i^*\}, B\Delta\{j,j^*\} \}\cap \B| < 2. \]
\end{enumerate}
\end{thm}
\begin{proof}
Two vertices of a polytope are adjacent exactly if there exists a
linear functional which takes equal, maximal values on the vertices,
and smaller values on all other vertices of the polytope.  Such a
functional in the first case is
\[ L=\sum_{k\in (A\cap B\cap I)} x_k - \sum_{k^*\in (A\cap B\cap I^*)}
x_i, \] by inspection. In the second case, we have the condition
that not both of the two sets listed are bases.  Suppose without
loss of generality that $C=B\Delta \{i,i^*\}$ is not a basis.
Then
\[ L=\sum_{k\in (C\cap I)} x_k - \sum_{k^*\in (C\cap I^*)} x_k +
\sum_{k\in (A \cap B \cap I)} x_k - \sum_{k^*\in (A \cap B \cap
I^*)} x_k
\] is clearly a required functional.
\end{proof}

We can show that in the Lagrangian case, the faces of the matroid polytope
can
only be the following:
\begin{description}
\item[sSquare] A square with short edges
(see, for example, the face of Figure 1 whose vertices all contain $2^*$).
\item[lSquare] A square with long edges  (an example follows).
\item[nsRect] A rectangle with two short and two long edges (see Figure
1).
\item[iTri] An isosceles triangle with two short and one long edge (see
Figure
1).
\item[eqTri] An equilateral triangle with long edges
(any face of Figure 2, defined later).
\end{description}
The matroid polytope of the representation
$$ \left( \begin{array}{cccc|cccc}
        0 & 1 & 0 &0& 1 & 0 & 0 &0\\
        1 & 0 & 0 &0& 0 & 1 & 0 &0\\
        0 & 0 & 0 &1& 0 & 0 & 1 &0 \\
        0&0&1&0& 0 & 0 & 0 & 1
\end{array} \right) $$
can be easily seen to be a square with long edges, embedded in 4-space.

Again, for details, see \cite{car}; these results all follow
without great difficulty from results in \cite{bgw} and
\cite{polytopes}.

\section{Oriented Lagrangian Matroids}

Let $M = (J,*,\B)$ be a Lagrangian matroid. Let $F \in \B$ be a
basis of $M$; we shall call it the \emph{fundamental basis}. The
fundamental basis defines a partial ordering on $J_n$ and thus on
$\B$ in the following way. The \emph{height} $h(A)$ of $A\in J_n$
relative to $F$ is defined to be $h(A) =n- | A \cap F |$. Although
$h(A)$ is dependent on $F$, we haven chosen not to emphasise this
fact in our notation. Notice that if we choose an admissible
ordering in which the smallest $n$ elements are those of $F$, then
$F$ is clearly the unique minimal basis in this ordering and this
admissible ordering never disagrees with the partial ordering; if
something is higher in the admissible ordering, then it is at
least as far from $F$ in height. Let $\Delta$ signify the
symmetric difference operator, $A\Delta B = (A\cup
B)\smallsetminus (A\cap B)$. Notice that a 2--dimensional face of
the hypercube has vertices of the form $\{ A, A \Delta \{i,i^*\},
A \Delta \{j,j^*\},A \Delta \{i,i^*,j,j^*\}\}$, and recall that
such faces contain all edges of the Lagrangian matroid polytope.
Then we see that such a face has one vertex closest to $F$, one
vertex furthest from $F$, and the other two vertices at the same
distance from $F$ with respect to height. We call an edge of the
matroid polytope between these last two vertices a
\emph{horizontal long edge (relative to $F$)}. The other possible
type of long edge is called a \emph{vertical long edge}. To
illustrate, we show Figure 2, where the fundamental basis
$1^*2^*3^*$ is shown as a larger dot than the other bases. In this
polytope, those edges incident with the fundamental basis are
vertical long edges, and the other edges are all horizontal long
edges. It can be shown that this is a symplectic matroid; however,
it is not representable in the sense used in this paper other than
over a field of characteristic 2. Over $\mathbb{GF}(2)$, it can be
represented as $$\left( \begin{array}{ccc|ccc} 0 & 1 & 1 & 1 & 0 &
0\\ 1 & 0 & 1 & 0 & 1 & 0\\ 1 & 1 & 0 & 0 & 0 & 1
\end{array}\right). $$

\begin{figure} \caption{}

\unitlength 1mm
\begin{picture}(80,80)(-20,0)
{\tiny
\put(51,9){\makebox(0,0)[tl]{$(1,1,-1)$}}
\put(9,51){\makebox(0,0)[br]{$(1,-1,1)$}}
\put(29,25){\makebox(0,0)[r]{$(-1,-1,-1)$}}
\put(71,66){$(-1,1,1)$}
}
\qbezier[20](10,50)(30,50)(50,50)
\qbezier[20](50,10)(50,30)(50,50)
\qbezier[20](10,10)(30,10)(50,10)
\qbezier[20](10,10)(10,30)(10,50)
\qbezier[20](30,65)(50,65)(70,65)
\qbezier[20](70,25)(70,45)(70,65)
\qbezier[20](30,25)(50,25)(70,25)
\qbezier[20](30,25)(30,45)(30,65)
\qbezier[15](50,50)(60,57.5)(70,65)
\qbezier[15](10,50)(20,57.5)(30,65)
\qbezier[15](10,10)(20,17.5)(30,25)
\qbezier[15](50,10)(60,17.5)(70,25)
\thicklines
\put(30,25){\line(1,1){1.5}}\put(33.5,28.5){\line(1,1){36.5}}
\put(30,25){\line(4,-3){20}}
\put(30,25){\line(-4,5){20}}
\put(10,50){\line(1,-1){40}}
\put(10,50){\line(4,1){60}}
\qbezier(50,10)(60,37.5)(70,65)
\put(50,10){\circle*{1.49}}
\put(10,50){\circle*{1.49}}
\put(30,25){\circle*{2.5}}
\put(70,65){\circle*{1.49}}
\thinlines
\end{picture}

\end{figure}

Let $s_F: \B \rightarrow \{ +1, -1 \} $ be a map from the
collection of bases of $M$ onto the two signs. We say that
$s_F(A)$ for $A \in \B$ is the \emph{sign} of the basis $A$
relative to the fundamental basis $F$. Then we say that $s_F$ is a
\emph{relative sign function} on $M$ with respect to $F$  if the
following four axioms are satisfied (in addition to the Maximality
Property). We shall say that $s_F$ defines an \emph{orientation}
on $M$. Abusing the terminology, we shall also call $s_F$ an
\emph{orientation of $M$ relative to $F$}, although this usage
does not agree well with the conventions of the theory of oriented
(ordinary) matroids.

\begin{ax} \label{axa}
If there is a horizontal long edge (relative to $F$) between
$A,B \in \B$ then $s_F(A)=s_F(B)$.
\end{ax}
\begin{ax} \label{axb}
If there is a vertical long edge (relative to $F$) between
$A,B \in \B$ then $s_F(A)=-s_F(B)$.
\end{ax}
\begin{ax} \label{axc}
In a square 2--dimensional face with short edges, (type \textbf{sSquare}),
if
three bases have one sign relative to $F$
and the fourth
the other, that fourth basis must be the highest or lowest basis with
respect to $F$ in the face.
\end{ax}
\begin{ax} \label{axd}
$s_F(F)=+1$.
\end{ax}

Notice that each of these axioms may be checked one
$2$-dimensional face of the hypercube at a time, as every edge
lies in a such a face and the short-sided squares of Axiom
\ref{axc} are such faces. We now discuss change of fundamental
basis. Assume $s_F$ is an orientation of $M$ relative to $F$.

\begin{defn} \label{exchdef}
For $A, G \in \B$ we define $s(G,A) = s_F (A) \cdot s_F (G) \cdot
(-1)^{|G \smallsetminus (F \cup A)|}$.
\end{defn}

\begin{thm} \label{exchworks}
The function $s: {\cal B} \times {\cal B} \rightarrow \{\,+1, -1\,\}$
defined above satisfies
\begin{enumerate}
\item For all $G \in {\cal B}$, $s_G(\,.\,)=s(G,\,.\,)$  is an orientation
of
$M$ relative to the fundamental basis $G$;
that is, it satisfies Axioms \ref{axa}--\ref{axd}.
\item For all $G,H,A\in {\cal B}$, we have
     $$ s(G,A)=s(H,A)s(H,G)(-1)^{|G\smallsetminus(H\cup A)|}.$$
\end{enumerate}
\end{thm}

This is easily proven from the axioms and definition.
We thus extend $s_F$ to a
function $s: J_n \times J_n \rightarrow \{\,+1, -1, 0\,\}$
by means of the definition above and $s=0$ whenever
either of its arguments falls outside $\B$.
We say that two orientations $s_F, s_G$ of $M$ relative
to fundamental bases $F, G$ are {\em equivalent\/} whenever
$s_F, s_G$ extend to the same $s$.

\begin{defn} \label{def:new}
We say that the pair\/ $(M, s)$, with\/ $ M \subseteq J_n$ and $s$ is a
function
$s: J_n \times J_n \rightarrow \{+1, -1, 0\}$,
 is an \emph{oriented Lagrangian matroid} whenever\/
$M$ is a Lagrangian matroid,
$s=0$ if one of its arguments is not a basis, and $s$ satisfies
the conclusions of Theorem~\ref{exchworks}.
We also say that $s$ is an orientation of $M$.
\end{defn}

Equivalently, an oriented Lagrangian matroid may be
regarded as an equivalence class of orientations relative to fundamental
bases.
We shall often suppress $s$ and refer to
$M$ as an oriented Lagrangian matroid.

\begin{defn}
An \emph{even} symplectic matroid is one in which the
number of starred elements in all bases have the same
parity (and so likewise for unstarred).
\end{defn}

\begin{cor}
Any even Lagrangian matroid is orientable, and there is only one possible
orientation.
\end{cor}

\begin{proof}
Choose a fundamental basis. Orient it with positive sign, and
where $h$ is the height function relative to this basis, assign
each other basis $A$ the sign $(-1)^{h(A)/2}$. Note that this is
$\pm 1$, as $h$ is always even. Axiom \ref{axa} is satisfied,
since horizontal edges connect only bases of the same height;
Axiom \ref{axb} is satisfied as vertical edges connect bases of
height differential $2$. There are no short edges, so Axiom
\ref{axc} is satisfied, and Axiom \ref{axd} was specifically
stated above, so this defines an orientation relative to the
fundamental basis chosen. It is immediate from consideration of
Axioms \ref{axa} and \ref{axb} that no other orientation is
possible.
\end{proof}

Thus, the theory of symplectic orientations of even Lagrangian
matroids is rather uninteresting. Notice that symmetric matrices
over $\R$ which give rise to even Lagrangian matroids are very
special. Indeed, it follows from the later result in this paper
(Corollary~\ref{cor:last}) that if $A$ is a real symmetric matrix
and all diagonal minors of $A$ of odd dimension equal to zero,
then  the corresponding quadratic form $Q$ has signature $0$ and
allows a hyperbolic basis, that is, $Q$  can be transformed to the
form $$ x_1x_2 + x_3x_4 + \cdots + x_{2k-1}x_{2k}.$$

 Even Lagrangian matroids are in fact just
Lagrangian orthogonal matroids (see \cite{orth}), and a more
interesting orientation scheme for Lagrangian orthogonal matroids
is discussed in \cite{boo2}, which is a development of the
orientation scheme given in \cite{w3} with the extra consideration
of what we have called fundamental basis exchange.

\section{Index}

\begin{defn}
Given a Lagrangian matroid $M$ of rank $n$, its matroid polytope $P$ and a
fundamental basis $F$, we define an \emph{increasing path} from a basis
$\alpha$ to a basis $\omega$ as a list of vertices $\alpha=v_0, v_1,
\ldots,
v_s=\omega$ and edges $e_1,e_2,\ldots,e_s$
of the matroid polytope such that
\begin{itemize}
 \item[{\rm (1)}] For $1\leqslant i\leqslant s$, the endpoints of $e_i$ are $v_{i-1}$ and
$v_i$; and
 \item[{\rm (2)}] There is admissible ordering $\prec$ on $J$
 such that $F$ is the $\prec$-minimal basis of $M$, and,  for all $1\leqslant i\leqslant s$,
  $A_{i-1} \prec A_i$
 for the bases $A_{i-1}$ and $A_i$ corresponding to the vertices $v_{i-1}$ and
 $v_i$.
\end{itemize}
It follows from \cite[Lemma~3.4]{property} that condition (2) is
equivalent to the following reformulation in terms of linear
functionals.
\begin{itemize}
 \item[{\rm (2*)}] There is some linear functional $L$ on $\R^n$ which is compatible
with
the height function $h$
(that is, $h(A)>h(B)$ implies that $L(A) >L(B)$ for any bases $A$ and $B$
of
$M$), and which satisfies $L(v_i)>L(v_{i-1})$
for $1\le i\le s$.
\end{itemize}
We say that the path $(\{v_i\},\{e_i\})$ is increasing with respect to
$L$.
\end{defn}

We can find an $L$ on $R^n$ which agrees with any admissible
ordering, and thus with a corresponding height function. Suppose
without loss of generality that the ordering is $n\succ n-1 \succ
\cdots\succ 2 \succ 1 \succ 1^* \succ \cdots \succ n^*$. Then we
may take $L(x_1 e_1 + \cdots + x_n e_n) = \sum_1^n 3^i x_i$ on
$\R^n$; this agrees with the ordering and never takes the same
value on different vertices. The existence of an increasing path
from a given non-maximal vertex $\alpha$ to some $\omega$ such
that $L(\omega)$ is the maximal value achieved by $L$ on the
matroid polytope is a standard fact of linear programming.

\begin{defn}
Given an oriented Lagrangian matroid $M$ and a fundamental basis $F$, we
define
the \emph{index} of $M$ relative to the $F$ as the number
of changes of sign in any increasing path from $F$ to a basis at maximal
height.
\end{defn}

\begin{thm}
For a given oriented Lagrangian matroid $M$, index relative to $F$
is well-defined; that is, the number of changes of sign in any
increasing path from the fundamental basis to any basis of maximal
height from it is the same. \label{th:index}
\end{thm}

We shall spend the remainder of the section proving this in stages.
We first prove that any two paths leading to the
same basis of maximal height have the same number of changes of sign. Fix
$F$,
the fundamental basis, until the end of the section.

\begin{lem}
If two such paths $e, f$ differ only as to the route taken through a
single
two-dimensional face $h$ of the matroid polytope,
then they have the same indices.
\end{lem}

\begin{proof} This can be checked, case by case, by considering the
restrictions that the axioms place on each of the possible face
types listed in the previous section. We exhibit the first the
case of the rectangle \textbf{nsRect}. Notice that short edges
connect two bases whose heights differ by 1; thus, the four
vertices of the face take at least two different heights. Let  a
vertex of the face of minimal height correspond to the basis $A$.
Then the other vertices are
$$ B=A\Delta \{i,i^*,j,j^*\}, \qquad
C=A\Delta \{k,k^*\}, \qquad D=A\Delta\{i,i^*,j,j^*,k,k^*\}$$
for
some $i,j,k \in I$ satisfying $|\{i,i^*,j,j^*,k,k^*\}|=6$. The
short edges are $AC$ and $BD$, the long edges $AB$ and $CD$.
Suppose there are two vertices of the same height. Then both long
exchanges are horizontal, and $s(F,A)=s(F,B)$, $s(F,C)=s(F,D)$
where $F$ is the fundamental basis, by Axiom 2. If $s(F,A)=s(F,C)$
also, the case becomes trivial, so assume otherwise. As paths must
be non-decreasing in height, a path through this face travels
exactly one short edge if it moves from one of $A,B$ to one of
$C,D$ and no short edges otherwise. As sign changes take place
exactly when traversing short edges, the case follows.

If there are not two vertices at the same height, then both long edges are
vertical, and thus induce sign changes.
Now $h(A)<h(C)<h(B)<h(D)$. The only journey for which there is a choice of
two
increasing paths is $A\rightarrow D$, which
may take place as $AC+CD$ or $AB+BD$; if there is no choice of increasing
paths, the case is trivial.
Axiom 3 gives $s(F,A)=-s(F,B)$ and $s(F,C)=-s(F,D)$. If $s(F,A)=s(F,C)$
then we
obtain $s(F,B)=s(F,D)$ and both paths
contain one sign change. Alternatively, if $s(F,A)=-s(F,C)$ then we obtain
$s(F,B)=-s(F,D)$ and both paths
contain two sign changes. This establishes the result in this case.

We now exhibit the case of the triangle \textbf{iTri}. Let $A, B$ be the
bases
at the ends of the long edge and $C$ the
third basis.  If the long edge is vertical, with say $A$ lowest, the only
possible paths are $AB$ and $ACB$. Since the edge
is vertical, $B$ is of the opposite sign to $A$, and so whichever sign $C$
takes, either path has one sign change.

Suppose now $AB$ is horizontal. Thus $A$ and $B$ take the same sign, and
since
any pair of paths must either both
enter or both leave at $C$ (depending on whether it is higher or lower
than
$AB$), the result is clear.

Similar reasoning proves the lemma for the other types of two-dimensional
faces.
\end{proof}

The next two lemmata are proven for general convex polytopes. An
\emph{admissible linear functional} is one which
is not constant along any edge of the convex polytope; we have exhibited
such a
functional which agrees with a given
admissible ordering earlier.

\begin{lem}
If $\alpha$ is a vertex of the convex polytope $P$ (not
necessarily a matroid polytope), $L$ an admissible linear
functional, and $e$ and $f$ edges of $P$ which are incident to
$\alpha$ and increasing for $L$ when we head away from $\alpha$,
then there exists a path of two-dimensional faces {\rm (} so that
successive faces in the path have an edge in common\/{\rm )}
incident to $\alpha$ connecting $e$ to $f$ and all on the
$L$-positive side of $\alpha$.
\end{lem}

\begin{proof}
First cut $P$ by the hyperplane $H$ which is constant for $L$ and
passes through $\alpha$, letting $P'$ be the resulting polytope
containing $e$ and $f$.  Now let $P''$ be the vertex figure of
$P'$ at $\alpha$. Now $P''$ has a facet $H''$ contained in $H$,
and vertices $e''$ and $f''$ corresponding to the original edges
$e$ and $f$.  Now all we want is a path of edges from $e''$ to $f''$
which misses $H''$;  those edges will then correspond to the
desired two-dimensional faces in $P$.
This now follows easily from Theorem 15.4 in \cite{bron}.
\end{proof}

\begin{lem}
Any two paths
$e=(\{v_i\},\{e_i\})$ , $f=(\{w_i\},\{f_i\})$ increasing relative
to the same $L$ from $v_0=w_0=\alpha$ to $v_s=w_t=\omega$
with $L$ maximal in the polytope at $\omega$
may be transformed into each other in steps, with
each step involving only a transformation within
a single 2--dimensional face of the polytope
(a \emph{transformation facewise}),
and with the path at each step still an increasing path relative to $L$.
\end{lem}

\begin{proof}
Set $P(\alpha)$ to be the assertion that the theorem holds for some
particular
$\alpha$ as above, and now attempt to prove $P(\alpha)$ for every vertex
in the
polytope (which is exactly the theorem) by an induction we shall call
$(*)$. We
must prove:
\begin{description}
\item[Basis of (*)] $P(\omega)$ holds.
\item[Inductive Step of (*)]
Given that $P(\gamma)$ holds for all vertices $\gamma $ with
$$L(\gamma) > L(\alpha),
$$ we may deduce $P(\alpha)$.
\end{description}
The basis is clearly trivial.

The previous lemma tells us that we can find a path of two-dimensional
faces
$h_1, \ldots, h_k$ of the polytope where each face $h_i$ has two edges
$g_{i-1}, g_i$ incident to $\alpha$, $g_0=e_1$ and $g_k=f_1$,
and all the $g_i$'s increase going away from $\alpha$.
Let $Q(k)$ be the assertion that $P(\alpha)$ holds when the number of
faces
above is exactly $k$. If we can prove $Q(k)$ for all k under the inductive
hypothesis of $(*)$, the theorem is complete. We thus assume $P(\gamma)$
for
all vertices $\gamma $ with
 $L(\gamma) > L(\alpha)$, and proceed by induction $(**)$ on $k$.
\begin{description}
\item[Basis of (**)] $Q(0)$ holds.
\item[Inductive Step of (**)] $Q(k-1)$ implies $Q(k)$.
\end{description}
Let $\beta$ be the point at the end of $g_k$.
If $k=0$, the two paths have the same first edge, and since $L(\beta) >
L(\alpha)$ we have $P(\beta)$; but the two paths differ only after
$\beta$, and
so $P(\alpha)$. Thus the basis of $(**)$ holds.

Take now $k>0$. If $\beta$ is the vertex of maximal $L$ in the face
$h_k$, then there is an increasing path segment from $\alpha$ through
$g_{k-1}$
around the edge of the face to $\beta$, which then gives
us an alternative path replacing $f$ through a transformation facewise
with $k$
replaced by $k-1$, completing the induction $(**)$.

If not, we can replace the path from $\beta$ to $\omega$ with one that
goes
through
the point of maximal $L$ in the face, using again $P(\beta)$, and now
continue
as above. This again gives a situation where $k$ is decreased by one after
a
transformation facewise, and so the induction is complete.
\end{proof}

These lemmata suffice to prove the theorem for all increasing paths
relative to
the same linear ordering finishing at the same basis of maximal height. If
we
can
 exhibit a path
strictly increasing with respect to height, that is increasing and without
horizontal edges, this is increasing relative to any
$L$ agreeing with height, and we will have the result for all increasing
paths
to the same basis of maximal height. Consider some $\omega$ of maximal
height,
and take $\alpha$
the fundamental basis.
It is now convenient to turn the situation `upside down', as $\omega$ may
not
be unique of maximal height, but $\alpha$ is certainly unique of minimal
height.
If we suppose for convenience that $\alpha= \{1^*, \ldots, n^*\} = (-1,
\ldots,
-1)$ as a point in the
space, and set a functional $L= -\sum_{i=1}^n x_i$ for a point with
co-ordinates $(x_i)$, then we see easily that $\alpha$
is maximal in $L$ and that $L$ is constant along horizontal edges. From
linear
programming, we know that from any non-maximal point
there extends a strictly increasing edge.
Thus we may find a never-horizontal path from $\omega$ to $\alpha$, which
is
the reverse of the path required.

It remains only to prove that the theorem holds for increasing paths
terminating at different bases of maximal height. Notice that any two
such bases may be connected by a path of horizontal long edges. It
thus suffices to show that two paths terminating at adjacent vertices
$\omega, \omega' $ of maximal height contain the same number of sign
changes. Now, we may replace both paths by never-horizontal paths $p,
p'$ respectively as above. Consider any admissible linear functional
$L$ agreeing with height (we have seen that such things exist). From
admissibility, without loss of generality $L(\omega) < L(\omega')$, as
no two vertices take equal values. Now consider the path $q$ to
$\omega' $ which consists of $p$ followed by the horizontal edge
between $\omega, \omega' $. Both $q$ and $p'$ are increasing relative
to $L$ and terminate at $\omega' $, and so they have the same number
of changes of sign. Thus $q$ and $p$ have the same number of changes
of sign also, as they differ only in the final horizontal edge. This
completes the proof.
\medskip

Notice that the above argument allows us to formulate an
elementary version of Theorem~\ref{th:index}.  Let $M$ be a
Lagrangian matroid of rank $n$, $P$ its matroid polytope and $F$ a
fundamental basis in $M$. Let $\alpha$ be the vertex of $P$
corresponding to the basis $F$. Then $\alpha$ can be connected
with any vertex $\omega$ of $P$ by a path $\alpha=v_0, v_1,
\ldots, v_s=\omega$ of vertices of the matroid polytope such that:
\begin{itemize}
 \item[{\rm (1)}] For $1\le i\le s$,  $v_{i-1}$ and
$v_i$ are adjacent in $P$; and
 \item[{\rm (2)}] The height is strictly increasing along the
 path:
 $h(v_{i-1}) < h(v_i)$ for all $1\le i\le s$.
\end{itemize}

We say that  $\alpha=v_0, v_1, \ldots, v_s=\omega$ is {\em
height-increasing} path.

\begin{thm}
For a given oriented Lagrangian matroid $M$, the number of changes
of sign in any height-increasing path from the fundamental basis
to any basis of maximal height from it is the same.
\label{th:height-index}
\end{thm}

We can take Theorem~\ref{th:height-index} for a more elementary
definition of index of the  oriented Lagrangian matroid $M$ with
respect to a fundamental basis.

\section{\relax Lagrangian Representations of\\
Oriented Lagrangian Matroids}

Recall that some Lagrangian matroids have symplectic representations
consisting
of two $n \times n$ matrices $A$
and $B$ indexed by $I,I^*$ respectively such that $AB^t$ is symmetric,
where a
set of $n$ columns of $(A,B)$ corresponds to a basis if and only if its
determinant
as a minor is non-zero. Such a representation also defines a Lagrangian
oriented matroid, as we shall now show.

Suppose $M$ is the Lagrangian matroid represented by $(A,B)$. Let
$F$ be a basis of $M$, the fundamental basis. Put all of the
columns indexed by elements of $F$ onto the right-hand side by
swapping $i$ with $i^*$ as required and multiplying the column
swapped into the left-hand side by $-1$. When all the columns of
$F$ are on the right, perform row operations to reduce the right
hand $n$ columns to the identity (possible as $F$ is a basis). Now
we have the form $(C_F, I_n)$ where $C_F$ is symmetric and
completely determined by $A, B$ and $F$. Observe that $K \in J_n$,
a set of $n$ column labels, is a basis if and only if the square
diagonal minor indexed by the columns of $K'$ is non-zero, where
$G'$ denotes $G \smallsetminus F$ for all $G\in J_n$. We define
$s(F,G)$ for $G \in \B$ to be the sign of the minor indexed by
$G'$.

We must now show that this definition satisfies the axioms for
$s(F,\,.\,)$ given in the previous section, and that choosing a
different fundamental basis gives results in accordance with
Definition \ref{exchdef}; that is, that Definition \ref{def:new}
holds. We first check agreement with Definition \ref{exchdef}, in the
form
\[ s(G,K) = s(F,K)\cdot s(F,G)\cdot(-1)^{|G\smallsetminus(F\cup K)|}.\]
 It is enough to prove that this is true in cases where $G=F \Delta
\{i, i^* \}$, as then the case where there are two exchanges (that is,
$G=F \Delta \{i, i^*, j, j^* \}$ with $i \ne j, j^*$) follows from the
continuity of the determinant function, and any path through a
Lagrangian matroid can be built up through these two types of
exchange. By this we mean that if both intervening sets $ F \Delta
\{i, i^* \}$, $ F \Delta \{j, j^* \}$ are non-bases, then we may make
small changes in entries of $(A,B)$ in order to make one or both
bases, and consider limits as these changes go to zero.

We consider two subcases, firstly the one where $K' \cap \{i,i^*\}
\neq \emptyset$. Now, it must be possible to reduce this $i$-th
column of $C_F$ to $e_i$, as otherwise $G$ is not a basis. In
doing so, we reduce what was $K'$ to one identity column and a
minor corresponding to the columns of $K \smallsetminus G$. Now
simple linear algebra tells us that the only change of the
determinant of this minor undergone is division by the value in
the $i$-th row and column, which has the sign $s(F,G)$. Hence $$
s(G,K) = s(F,G) \cdot s(F,K) $$ which is what we require as $G
\smallsetminus (F \cup K) = \emptyset$ in this case.

If $K' \cap \{i,i^*\} = \emptyset$, consider just the minor $K'
\cup {i}$, where without loss of generality $i \in G
\smallsetminus F$, and the corresponding identity minor in the
same position on the right hand side. Upon swapping columns $i$
with $i^*$ and multiplying by column $i^*$ by $-1$, the left hand
minor is $K'$ with $-1$ times the appropriate identity column and
row added, and so has determinant $-1$ times the determinant of
$K'$. The right hand minor is an identity matrix other than column
$i$, and so has determinant equal to the entry from the $i$-th row
of the $i$-th column which is $s(F,G)$; this is non-zero, as $G$
is a basis and this right-hand side represents $G$. When the right
hand side is reduced to the identity, the left hand side thus
becomes the minor of $K$ relative to $G$, and so will have a
determinant of sign $$ s(G,K) = s(F,G) \cdot s(F,K) \cdot (-1) $$
which is what we require as $G \smallsetminus (F \cup K) = \{ i
\}$. So Definition \ref{exchdef} is respected.

Axiom \ref{axd} is trivially satisfied, as the sign of the empty minor
is positive by definition.

For axiom \ref{axa}, consider a horizontal long edge from $K$ to
$L$. Notice that $| F \smallsetminus (K \cup L) |$ is of the
opposite parity to $| F \smallsetminus (K \cup K)
|=|F\smallsetminus K|$. Indeed, this holds if and only if
$|F\cap(K\cup L)|$ is of the opposite parity to $|F\cap K|$.  Now,
\[ F\cap(K\cup L) = (F\cap K)\cup(F\cap L) = (F\cap K)\sqcup
(F\cap(L\smallsetminus K)).\] Thus, we need to show that
$F\cap(L\smallsetminus K)$ is of odd parity. Also, $|(F\cap
L)|=|(F\cap K)|$ as the edge is horizontal, and (similarly)
\[ F\cap(K\cup L) = (F\cap L)\sqcup(F\cap(K\smallsetminus L)). \]
So $|F\cap(K\smallsetminus L)|=|F\cap(L\smallsetminus K)|$.
Since, for appropriate $i,j\in J$, we have $K \Delta L =
\{i,i^*,j,j^*\}$, we can write
\[ (F\cap(K\smallsetminus L))\sqcup (F\cap(L\smallsetminus K)) = F\cap (K\Delta
L)
 = F\cap \{i,i^*,j,j^*\} \]
which has size 2, since $F$ is maximal admissible.  This establishes
 that the set $ F \smallsetminus (K \cup K)$ has opposite parity to $ F \smallsetminus
 (K \cup L)$.

Now $$ s(F,K) = s(K,F) \cdot (-1)^{| F \smallsetminus (K \cup K)
|}$$ and $$ s(F,L) = s(K,L) \cdot s(K,F) \cdot (-1)^{| F
\smallsetminus (K \cup L) |}$$ and so $ s(F,K) = s(F,L) $ if and
only if $s(K,L)=-1$. But this holds, as the minor corresponding to
$s(K,L)$ where there is a long edge is of the form $$ \left(
\begin{array}{cc} a & b \\ b & 0 \end{array} \right) $$ which must
have negative sign. It must take this form because at least one of
the vertices of the $n$--cube adjacent (in the $n$-cube) to the
fundamental basis and lying within the same 2--dimensional face of
the $n$--cube as this long edge must be a non -basis, otherwise
the long edge would not exist; the zero corresponds to this
vertex. Axiom \ref{axb}, concerning a vertical long edge from $K$
to $L$, is similarly satisfied through the observance that here $|
F \smallsetminus (K \cup L) |$ must be of the same parity as $| F
\smallsetminus (K \cup K) |$, and the rest of the argument
similar.

It only remains to prove that axiom \ref{axc} holds.  Let $K$ be the
basis of minimal height in a short-sided square face; then the other
bases can be written as
\[ L_1 = (K\smallsetminus\{f_1\})\sqcup\{f_1^*\}, \quad
   L_2 = (K\smallsetminus\{f_2\})\sqcup\{f_2^*\}, \quad
   L_3 = (K\smallsetminus\{f_1,f_2\})\sqcup\{f_1^*,f_2^*\} \]

for some $f_1,f_2\in F$.  We show first that the relative signs of $L_3$
and $K$ and of $L_1$ and $L_2$ are unchanged when $K$ replaces $F$ as
fundamental base.  Observe first that
\[ K\smallsetminus (F\cup K) = K\smallsetminus (F\cup L_1)
= K\smallsetminus (F\cup L_2) = K\smallsetminus (F\cup L_3) =
\emptyset; \] now, from Definition \ref{exchdef} and the above, we
obtain \( s(K,L_i) = s(F,L_i) \cdot s(F,K) \) for $i=1,2,3$, and
from this \(s(K,L_1)\cdot s(K,L_2) = s(F,L_1)\cdot s(F,L_2)\) as
required.  Thus these relative signs are unchanged by change of
fundamental base; and to prove that the axiom holds, it suffices
to show that if the bases at the same height are of opposite sign,
so are the other two.

Now consider a minor corresponding to this face with $K$ as the
fundamental base where the bases of equal height take
different signs; it must take the form
$$ \left( \begin{array}{cc} a & b \\ b & -c \end{array} \right) $$
where $a$ and $c$ are non-zero and of the same sign. But then the
minor of the top basis, which is shown, is clearly negative, but the
fundamental basis is always positive. So axiom \ref{axc} holds. Thus
we have proven
\begin{thm}
A representation of a Lagrangian matroid $M$ defines also an
orientation of $M$, by the above procedure.
\end{thm}
We make the obvious definition: an $n \times 2n$ matrix $(A,B)$ with
$AB^t$ symmetric is called a \emph{representation} of the oriented
Lagrangian matroid it yields through the above procedure.

Since, in general, there is more than one way to orient a
Lagrangian matroid, we see that Schubert cells are partitioned by
orientation. It is clear that different orientations of the same
Lagrangian matroid lie in different connected components of the
thin Schubert cell; thus, Lagrangian orientations provide a
natural further refinement of the structure on a real Lagrangian
variety on a symplectic space.

Recall that in linear algebra, the index of a quadratic form is
often defined as the number of negative coefficients when it is
expressed as a sum of squares. The signature, more often referred
to, is the number of positive terms less the index, so that index
can be recovered from rank and signature. It is well known that
rank, signature and index are all invariant under change of basis
of a quadratic form. We now give the theorem that motivates the
naming of the matroid property `index':

\begin{thm}
The index of a (represented) oriented Lagrangian matroid relative to a
fundamental basis $F$,
as defined in section 2, is the same
as the index of the quadratic form $Q$ determined by the symmetric matrix
$C_F$
constructed at the beginning of this section. Also, the rank of $Q$ is the
same
as the maximal height relative to $F$ attained in the matroid.
\end{thm}
\begin{proof}
Notice first that upon obtaining the form $(C_F, I_n)$ as above, bases of
maximal height correspond to non-zero minors of $C_F$ of maximal
size. This
means exactly that rank of $Q$ is this height.
 We utilise two results from the standard theory of quadratic forms (see,
for example, \cite[pages 133--138 in chapter 17]{mat} for a
particularly clear treatment. Note, however, that Ayres defines
index as the number of positive terms, rather than negative). If
$P_i$ is the determinental minor consisting of the first $i$
columns and rows of a square, symmetric matrix $A$ then we can
rearrange by swapping of columns and of the corresponding rows so
that not both of $P_{i}$ and $P_{i-1}$ are zero for $1\leqslant i
< r$, and $P_r$ is non-singular, where $r= \hbox{rank } A$. A
matrix where this holds is called \emph{regularly arranged}.
Notice that putting  $A$ into a regular arrangement, and permuting
column labels when permuting columns (of A), does not alter the
oriented Lagrangian matroid represented by $(A, I_n)$. When a
matrix is in regular arrangement, Kronecker's method tells us that
the index is exactly the number of changes of sign in the sequence
$P_0, P_1, \ldots, P_r$, discounting any zeroes.

We may assume that $C_F$ is in a regular arrangement.
We choose an increasing path by taking the list of
bases
corresponding to the non-zero $P_i$; this is clearly increasing, and in
fact
has no horizontal edges. The theorem is immediate.
\end{proof}

\section{Oriented Lagrangian Matroid Polytopes}

Take a Lagrangian matroid polytope, and assign a direction to each short
edge,
allowing long edges to remain undirected. Then
\begin{defn}
A Lagrangian matroid polytope with oriented short edges is an
\emph{oriented
Lagrangian
matroid polytope} exactly
when every 2--dimensional face has the same number of short edges
directed
clockwise as
anti-clockwise.
\end{defn}

Assign  a height function relative to some vertex of the $n$--cube (height
is
defined relative to non-bases by simply dropping the requirement that $F$
is a
basis in the original definition of height). Then an \emph{inducing edge}
relative to this height
function is one which is either (long and) vertical, or which is (short
and)
directed downwards (notice all short edges connect bases of heights
differing by 1).

\begin{thm} \label{thm:oripol}
Take an oriented Lagrangian matroid polytope and choose a
fundamental basis $F$ and corresponding height function. Assign
signs $s(F,\,.\,)$ as follows: $s(F,F)$ is positive. Any two bases
connected by inducing edges have opposite signs. Any two bases
connected by non-inducing edges have the same signs. Now:
\begin{enumerate}
\item This procedure is contradiction-free; and
\item The result is an oriented Lagrangian matroid relative to $F$, and
different fundamental bases give equivalent orientations, so that an
oriented
Lagrangian matroid polytope defines uniquely an oriented Lagrangian
matroid.
\item Every oriented Lagrangian matroid is induced from an orientation of
its Lagrangian matroid polytope, so that oriented polytopes and
oriented Lagrangian  matroids are in one to one correspondence.
\end{enumerate}
\end{thm}

\begin{proof}
We first consider contradiction-freeness.
The result is trivial when the Lagrangian matroid polytope is
one-dimensional
or smaller. If not, it is enough to show
the procedure contradiction--free on 2--dimensional faces, which is a
simple
check.
For part two,
as observed when  axioms
\ref{axa}--\ref{axd}
were stated, it is enough to check axioms on these faces also, and this
too is
immediate. It remains only to show that
choosing a different fundamental basis gives results respecting Definition
\ref{exchdef}. The only possible points of concern
are around short edges, since long edges give no freedom of sign choice if
the
result is a oriented Lagrangian matroid, as we have shown it must
be. Notice
that where $A, B \in \B$ are the ends of a short edge, that is $A\Delta B
=
\{i, i^*\}$ for some $i$, Definition \ref{exchdef} yields $ s(A,B)
= - s(B,A) $ regardless of what signs we are given relative to
$F$. Fix some $F$, and take all signs as being calculated from the
signs relative to $F$ decided as above through Definition
\ref{exchdef}. Let us suppose the edge is directed  from $A$ to
$B$. We first show $s(A,B)=1$ regardless of signs relative to $F$.
Now, $$ s(A,B) = s(F,A) s(F,B) (-1)^{|A\smallsetminus (B\cup
F)|}.$$ But if $A$ is lower than $B$, then $|A\smallsetminus
(B\cup F)|=0$ and $s(F,A)=s(F,B)$, as the edge is non-inducing, so
$s(A,B)$ is positive as required. If $A$ is higher then $B$,
$|A\smallsetminus (B\cup F)|=1$ and $s(F,A)=-s(F,B)$, so $s(A,B)$
remains positive. Now consider some basis $G$ as fundamental. To
complete the proof of this point, we must show that $s(G,B) s(G,A)
= 1$ if and only if $A$ is lower than $B$ relative to $G$. But $$
1 = s(A,B) = s(G,A) s(G,B) (-1)^{|A\smallsetminus (B\cup G)|}$$
    from the above and Definition \ref{exchdef}. But notice, as before,
that $|A\smallsetminus (B\cup G)|$ is 0 when $A$ lies below $B$
relative to $G$ and 1 otherwise, which completes our proof.
Finally, for part three, observe once more that only short edges
need give us concern. From Definition \ref{exchdef} we obtain
$s(A,B) = -s(B,A)$ whenever $A, B$ are the ends of a short edge,
and so we direct the edge according to which of the two is
positive. This then obviously
 gives the correct oriented polytope.
\end{proof}

\begin{cor} \label{cor:last}
The index of an oriented Lagrangian matroid relative to some basis $F$ may
be
found by considering its
oriented polytope and counting inducing edges in any path from $F$ to a
basis
of maximal height in which
no edge is followed downwards relative to $F$.
\end{cor}

Notice that in an even Lagrangian matroid inducing edges are
vertical long edges. Hence the index of an even Lagrangian matroid
is with respect to the fundamental basis $F$ is exactly half of
the maximal height of bases with respect to $F$. For an even
Lagrangian  matroids represented by a symmetric matrix $A$ over
$\R$, this means that the index of the corresponding quadratic
form $Q$ is half its rank. Therefore $Q$ has signature $0$, hence
allows a hyperbolic basis, that is, can be transformed to the form
$$ x_1x_2 + x_3x_4 + \cdots + x_{2k-1}x_{2k}.$$

\begin{figure} \caption{}

\unitlength .8mm
\begin{picture}(80,80)(0,0)
{ \tiny
\put(51.5,49.5){\makebox(0,0)[tl]{$123$}}
}
\qbezier[20](10,10)(30,10)(50,10)
\qbezier[20](10,10)(10,30)(10,50)
\qbezier[20](30,65)(50,65)(70,65)
\qbezier[20](70,25)(70,45)(70,65)
\qbezier[2](30,25)(32,25)(34,25)\qbezier[7](56,25)(63,25)(70,25)
\qbezier[2](30,25)(30,27)(30,29)\qbezier[4](30,57)(30,61)(30,65)
\qbezier[15](10,50)(20,57.5)(30,65)
\qbezier[15](10,10)(20,17.5)(30,25)
\qbezier[15](50,10)(60,17.5)(70,25)
\thicklines
\put(10,50){\line(1,-1){40}}
\put(10,50){\line(4,1){60}}
\qbezier(50,10)(60,37.5)(70,65)
\put(10,50){\line(1,0){40}}
\put(50,10){\line(0,1){40}}
\put(50,50){\line(4,3){20}}
\put(30,50){\vector(1,0){0}}
\put(50,30){\vector(0,1){0}}
\put(60,57.5){\vector(-4,-3){0}}
\put(50,10){\circle*{1.49}}
\put(10,50){\circle*{1.49}}
\put(50,50){\circle*{1.49}}
\put(70,65){\circle*{1.49}}
\thinlines
\end{picture}
\begin{picture}(80,80)(0,0)
{ \tiny
\put(51.5,49.5){\makebox(0,0)[tl]{$123$}}
}
\qbezier[20](10,10)(30,10)(50,10)
\qbezier[20](10,10)(10,30)(10,50)
\qbezier[20](30,65)(50,65)(70,65)
\qbezier[20](70,25)(70,45)(70,65)
\qbezier[2](30,25)(32,25)(34,25)\qbezier[7](56,25)(63,25)(70,25)
\qbezier[2](30,25)(30,27)(30,29)\qbezier[4](30,57)(30,61)(30,65)
\qbezier[15](10,50)(20,57.5)(30,65)
\qbezier[15](10,10)(20,17.5)(30,25)
\qbezier[15](50,10)(60,17.5)(70,25)
\thicklines
\put(10,50){\line(1,-1){40}}
\put(10,50){\line(4,1){60}}
\qbezier(50,10)(60,37.5)(70,65)
\put(10,50){\line(1,0){40}}
\put(50,10){\line(0,1){40}}
\put(50,50){\line(4,3){20}}
\put(30,50){\vector(-1,0){0}}
\put(50,30){\vector(0,-1){0}}
\put(60,57.5){\vector(4,3){0}}
\put(50,10){\circle*{1.49}}
\put(10,50){\circle*{1.49}}
\put(50,50){\circle*{1.49}}
\put(70,65){\circle*{1.49}}
\thinlines
\end{picture}

\end{figure}
We present the examples shown in figure 3. Both polytopes
correspond to the same unoriented Lagrangian matroid, with bases
$1^*23, 12^*3, 123^*, 123$. If we consider signs relative to
$123$, all other signs are negative in the left-hand example and
all signs are positive in the right-hand example. Both these
oriented Lagrangian matroids are representable, by $$
\bordermatrix{ & 1 & 2 & 3 & 1^* & 2^* & 3^* \cr & 1 & 0 & 0 & 1 &
1 & 1 \cr & 0 & 1 & 0 & 1 & 1 & 1 \cr & 0 & 0 & 1 & 1 & 1 & 1 }
\sim \bordermatrix{ & 1^* & 2^* & 3^* & 1 & 2 & 3  \cr & -1 & -1 &
-1 & 1 & 0 & 0\cr & -1 & -1 & -1 & 0 & 1 & 0\cr & -1 & -1 & -1 & 0
& 0 & 1 } $$ and $$ \bordermatrix{ & 1 & 2 & 3 & 1^* & 2^* & 3^*
\cr & 1 & 0 & 0 & -1 & -1 & -1 \cr & 0 & 1 & 0 & -1 & -1 & -1 \cr
& 0 & 0 & 1 & -1 & -1 & -1 } \sim \bordermatrix{ & 1^* & 2^* & 3^*
& 1 & 2 & 3  \cr & 1 & 1 & 1 & 1 & 0 & 0\cr & 1 & 1 & 1 & 0 & 1 &
0\cr & 1 & 1 & 1 & 0 & 0 & 1 } $$
 respectively, where $\sim $ means that these are the
same representations after permutations of columns and column
labels as defined in section 3. All the indices of both examples
are 1 except for the index of the right-hand example relative to
$123$, which is 0. Finally, it is easy to see (from the definition
of oriented polytope and the theorem above) that these are the
only possible orientations of this polytope, and hence of this
Lagrangian matroid.

Finally, we consider isomorphisms of oriented Lagrangian matroids.
We say that two symplectic matroids on $J = I \cup I^* $ with
basis collections $\B, \B'$ are isomorphic when there exists some
admissible permutation $\sigma $ such that $$ \B' = \sigma \B = \{
\sigma A | A \in \B \}.$$ We say that two oriented Lagrangian
matroids, using the above notation with orientations $s, s'$ are
isomorphic when, in addition, $$ s'({\sigma A}, \sigma B) = s(A,
B) \hbox{ for all } A, B \in J_n.$$ Observe that two oriented
Lagrangian matroids are isomorphic exactly when their oriented
polytopes are the same up to some admissible permutation of
vertices of the $n$-cube, where an edge directed from $A$ to $B$
becomes an edge directed from $\sigma A$ to $\sigma B$. Recall
that the group of all admissible permutations is exactly the
symmetry group of the $n$-cube. Thus isomorphism is a very natural
geometric concept. We now see that the two Lagrangian matroid
polytopes of figure 3 are (trivially) isomorphic as unoriented
polytopes, but non-isomorphic as oriented polytopes.

\paragraph{Acknowledgement.} The authors wish to thank the referee
for comments which helped to improve the exposition in this paper.

\vfill

\noindent
{\bf Authors' addresses}

\bigskip
\noindent
Richard F. Booth, 

\bigskip
\noindent
Alexandre V. Borovik, Department of Mathematics, UMIST, PO Box 88,
Manchester M60 1QD, United Kingdom;
{ \tt
alexandre.borovik@umist.ac.uk
}

\bigskip
\noindent
Israel M. Gelfand, Department of Mathematics, Rutgers
University, New Brunswick NJ 08903, USA; {\tt igelfand@math.rutgers.edu}

\bigskip
\noindent
Neil White, Department of Mathematics, P.O. Box 118105, University of
Florida,
Gainesville FL 32611-8105, USA;
{ \tt
white@math.ufl.edu
}

\end{document}